\documentclass{article}%
\usepackage{amsfonts}
\usepackage{amsmath}%
\usepackage{amssymb}%
\usepackage{graphicx}
\newtheorem{theorem}{Theorem}

\newenvironment{proof}[1][Proof]{\noindent\textbf{#1.} }{\ \rule{0.5em}{0.5em}}

\begin{document}

\title{One Dimensional Metrical Geometry}
\author{N.\ J. Wildberger\\School of Mathematics and Statistics\\UNSW Sydney 2052 Australia\\web pages: http://web.maths.unsw.edu.au/\symbol{126}norman/}
\maketitle
\date{}

\begin{abstract}
One dimensional metrical geometry may be developed in either an affine or
projective setting over a general field using only algebraic ideas and
quadratic forms. Some basic results of universal geometry are already present
in this situation, such as the Triple quad formula, the Triple spread formulas
and the Spread polynomials, which are universal analogs of the Chebyshev
polynomials of the first kind. Chromogeometry appears here, and the related
metrical and algebraic properties of the projective line are brought to the fore.

\end{abstract}

\section{Introduction}

This paper introduces the metrical geometry of \textit{one dimensional space},
a surprisingly rich but largely unexplored subject. This theory is a natural
chapter of \textit{universal geometry}, introduced in \cite{Wild}, where all
fields not of characteristic two are treated uniformly.

Over the `real numbers' it is usual to consider the two models of one
dimensional geometry to be the \textit{line} and the \textit{circle}. Either
way one must struggle with the proper definition of the `continuum'. In this
paper I wish to promote instead the idea that the two fundamental one
dimensional objects are the \textit{affine line} and the \textit{projective
line}, and that a purely algebraic approach to metrical structure simplifies
the logical coherence of the subject.

The affine line has a unique metrical structure (not a metric in the usual
sense!) derived from the notion of \textit{quadrance }between two points. This
is the square of the distance in the usual `real number' sense, but its
algebraic nature allows it to be defined over a general field.

The projective line has a family of metrical structures determined by
\textit{forms,} with \textit{projective quadrance} playing the role of
quadrance. Projective quadrance between projective points is equivalent to
\textit{spread} between lines in the plane, and many of the results of this
paper are basic for \textit{rational trigonometry}, as described in
\cite{Wild}.

The formulas developed here also anticipate results of planar Euclidean
geometry, including Heron's or Archimedes' formula on the area of a triangle
and Brahmagupta's formula for the area of a cyclic quadrilateral. They extend
to affine and projective metrical geometries associated to arbitrary quadratic
forms in higher dimensions, as described in \cite{Wild3}, and so become
important for elliptic and hyperbolic geometries.

Of particular interest in the one dimensional situation are three forms
called\textit{\ blue}, \textit{red }and \textit{green, }as these lead to
\textit{chromogeometry, }an entirely new framework for our understanding of
planar metrical geometry, as described in \cite{Wild2}. The blue form
corresponds to Euclidean geometry, while the red and green forms are
hyperbolic, or relativistic, analogs. Each yields an \textit{algebraic
structure }to the projective line, and these three geometries interact in a
pleasant way even in the one dimensional situation.

The metrical and algebraic properties of the projective line introduced here
should be of interest to algebraic geometers, as the natural extension of
these ideas to higher dimensions allows \textit{metrical analysis} of
projective curves and varieties---still completely in the algebraic context.

\section{The affine line}

\subsection{Quadrance and Triple quad formula}

Fix a field, whose elements are called \textbf{numbers}. A \textbf{point}
$A\equiv\left[  x\right]  $ is a number enclosed in square brackets.

The \textbf{quadrance} $Q\left(  A_{1},A_{2}\right)  $ between the points
$A_{1}\equiv\left[  x_{1}\right]  $ and $A_{2}\equiv\left[  x_{2}\right]  $ is
the number
\[
Q\left(  A_{1},A_{2}\right)  \equiv\left(  x_{2}-x_{1}\right)  ^{2}.
\]

Clearly
\[
Q\left(  A_{1},A_{2}\right)  =Q\left(  A_{2},A_{1}\right)  ,
\]
and $Q\left(  A_{1},A_{2}\right)  $ is zero precisely when $A_{1}=A_{2}.$

The next result is fundamental, with implications throughout universal geometry.

\begin{theorem}
[Triple quad formula]For points $A_{1},A_{2}$ and $A_{3},$ define $Q_{1}\equiv
Q\left(  A_{2},A_{3}\right)  $, $Q_{2}\equiv Q\left(  A_{1},A_{3}\right)  $
and $Q_{3}\equiv Q\left(  A_{1},A_{2}\right)  $. Then
\[
\left(  Q_{1}+Q_{2}+Q_{3}\right)  ^{2}=2\left(  Q_{1}^{2}+Q_{2}^{2}+Q_{3}%
^{2}\right)  .
\]

\end{theorem}

\begin{proof}
First verify the polynomial identity%
\[
\left(  Q_{1}+Q_{2}+Q_{3}\right)  ^{2}-2\left(  Q_{1}^{2}+Q_{2}^{2}+Q_{3}%
^{2}\right)  =4Q_{1}Q_{2}-\left(  Q_{1}+Q_{2}-Q_{3}\right)  ^{2}.
\]
Suppose that $A_{1}\equiv\left[  x_{1}\right]  $, $A_{2}\equiv\left[
x_{2}\right]  $ and $A_{3}\equiv\left[  x_{3}\right]  $. Then%
\begin{align*}
Q_{1}+Q_{2}-Q_{3}  & =\left(  x_{3}-x_{2}\right)  ^{2}+\left(  x_{3}%
-x_{1}\right)  ^{2}-\left(  x_{2}-x_{1}\right)  ^{2}\\
& =2\left(  x_{3}-x_{2}\right)  \left(  x_{3}-x_{1}\right)  .
\end{align*}
So
\[
\left(  Q_{1}+Q_{2}-Q_{3}\right)  ^{2}=4Q_{1}Q_{2}.
\]

\end{proof}

Motivated by the Triple quad formula, define \textbf{Archimedes' function
}$A\left(  a,b,c\right)  $ for numbers $a,b$ and $c$ by%
\[
A\left(  a,b,c\right)  \equiv\left(  a+b+c\right)  ^{2}-2\left(  a^{2}%
+b^{2}+c^{2}\right)  .
\]
A set $\left\{  a,b,c\right\}  $ is a\textbf{\ quad triple} precisely when
$A\left(  a,b,c\right)  =0.$ Note that $A\left(  a,b,c\right)  $ is a
symmetric function of $a,b$ and $c,$ and that%
\begin{align*}
&  A\left(  a,b,c\right)  =4ab-\left(  a+b-c\right)  ^{2}\\
&  =2\left(  ab+bc+ca\right)  -\left(  a^{2}+b^{2}+c^{2}\right) \\
&  =4\left(  ab+bc+ca\right)  -\left(  a+b+c\right)  ^{2}\\
&  =%
\begin{vmatrix}
2a & a+b-c\\
a+b-c & 2b
\end{vmatrix}
\\
&  =-%
\begin{vmatrix}
0 & a & b & 1\\
a & 0 & c & 1\\
b & c & 0 & 1\\
1 & 1 & 1 & 0
\end{vmatrix}
.
\end{align*}

\begin{theorem}
[Heron's identity]Suppose that $D_{1}=d_{1}^{2},$ $D_{2}=d_{2}^{2}$ and
$D_{3}=d_{3}^{2}$ for some numbers $d_{1},d_{2}$ and $d_{3}.$ Then%
\[
A\left(  D_{1},D_{2},D_{3}\right)  =\left(  d_{1}+d_{2}+d_{3}\right)  \left(
-d_{1}+d_{2}+d_{3}\right)  \left(  d_{1}-d_{2}+d_{3}\right)  \left(
d_{1}+d_{2}-d_{3}\right)  \allowbreak.
\]

\end{theorem}

\begin{proof}
Observe that%
\begin{align*}
& \left(  d_{1}+d_{2}+d_{3}\right)  \left(  -d_{1}+d_{2}+d_{3}\right)  \left(
d_{1}-d_{2}+d_{3}\right)  \left(  d_{1}+d_{2}-d_{3}\right) \\
& =\left(  \left(  d_{1}+d_{2}\right)  ^{2}-d_{3}^{2}\right)  \left(
d_{3}^{2}-\left(  d_{1}-d_{2}\right)  ^{2}\right) \\
& =\left(  \left(  d_{1}+d_{2}\right)  ^{2}+\left(  d_{1}-d_{2}\right)
^{2}\right)  D_{3}-\left(  d_{1}+d_{2}\right)  ^{2}\left(  d_{1}-d_{2}\right)
^{2}-D_{3}^{2}\\
& =2\left(  D_{1}+D_{2}\right)  D_{3}-\left(  D_{1}-D_{2}\right)  ^{2}%
-D_{3}^{2}\\
& =\left(  D_{1}+D_{2}+D_{3}\right)  ^{2}-2\left(  D_{1}^{2}+D_{2}^{2}%
+D_{3}^{2}\right)  .
\end{align*}

\end{proof}

By the classical \textit{Heron's formula}, known by Archimedes, the expression
$A\left(  d_{1}^{2},d_{2}^{2},d_{3}^{2}\right)  $ is sixteen times the square
of the area of a triangle with side lengths $d_{1},d_{2}$ and $d_{3}.$

\subsection{Quadruple quad formula}

The next two results extend the Triple quad formula to four points.

\begin{theorem}
[Two quad triples]Suppose that $\left\{  a,b,x\right\}  $ and $\left\{
c,d,x\right\}  $ are both quad triples. Then%
\[
\left(  \left(  a+b+c+d\right)  ^{2}-2\left(  a^{2}+b^{2}+c^{2}+d^{2}\right)
\right)  ^{2}=64abcd.
\]
Furthermore if $a+b\neq c+d$ then%
\[
x=\frac{\left(  a-b\right)  ^{2}-\left(  c-d\right)  ^{2}\allowbreak}{2\left(
a+b-c-d\right)  }.
\]

\end{theorem}

\begin{proof}
Suppose that $\left\{  a,b,x\right\}  $ and $\left\{  c,d,x\right\}  $ are
quad triples, so that%
\begin{align}
\left(  x-a-b\right)  ^{2} &  =4ab\label{TwoQuads}\\
\left(  x-c-d\right)  ^{2} &  =4cd.\label{TwoQuads2}%
\end{align}
Take the difference between these two equations to yield%
\[
2x\left(  a+b-c-d\right)  +\left(  c+d\right)  ^{2}-\left(  a+b\right)
^{2}=4cd-4ab.
\]
If $a+b\neq c+d$ then
\[
x=\frac{\left(  a-b\right)  ^{2}-\left(  c-d\right)  ^{2}\allowbreak}{2\left(
a+b-c-d\right)  }.
\]
Substitute this into (\ref{TwoQuads}) or (\ref{TwoQuads2}) to get%
\begin{equation}
\left(  \left(  a-b\right)  ^{2}-\left(  c-d\right)  ^{2}-2\left(
a+b-c-d\right)  \left(  a+b\right)  \right)  ^{2}-16ab\left(  a+b-c-d\right)
^{2}=0.\label{QuadrupleQuad}%
\end{equation}
The left hand side of this equation can be rearranged to get
\[
\left(  \left(  a+b+c+d\right)  ^{2}-2\left(  a^{2}+b^{2}+c^{2}+d^{2}\right)
\right)  ^{2}-64abcd.
\]

On the other hand if $a+b=c+d$ then by (\ref{TwoQuads}) we see that $4ab=4cd,
$ so that $\left(  a-b\right)  ^{2}=\left(  c-d\right)  ^{2},$ and then
(\ref{QuadrupleQuad}) reduces also to zero.
\end{proof}

The \textbf{Quadruple quad function} $Q\left(  a,b,c,d\right)  $ is defined
by
\[
Q\left(  a,b,c,d\right)  \equiv\left(  \left(  a+b+c+d\right)  ^{2}-2\left(
a^{2}+b^{2}+c^{2}+d^{2}\right)  \right)  ^{2}-64abcd.
\]

It is interesting that the expression
\[
\left(  a+b+c+d\right)  ^{2}-2\left(  a^{2}+b^{2}+c^{2}+d^{2}\right)
\]
appears in Descartes' circle theorem, a result motivated by a question of
Princess Elisabeth of Bohemia around 1640, and rediscovered by both Beecroft
and Soddy (see \cite{Pedoe}, \cite{Coxeter}).

\begin{theorem}
[Quadruple quad formula]For points $A_{1},A_{2},A_{3}$ and $A_{4},$ define the
quadrances $Q_{ij}\equiv Q\left(  A_{i},A_{j}\right)  $ for all $i,j=1,2,3$
and $4$. Then%
\[
Q\left(  Q_{12},Q_{23},Q_{34},Q_{14}\right)  =0.
\]
Furthermore%
\begin{align*}
Q_{13} &  =\frac{\left(  Q_{12}-Q_{23}\right)  ^{2}-\left(  Q_{34}%
-Q_{14}\right)  ^{2}}{2\left(  Q_{12}+Q_{23}-Q_{34}-Q_{14}\right)  }\\
Q_{24} &  =\frac{\left(  Q_{23}-Q_{34}\right)  ^{2}-\left(  Q_{12}%
-Q_{14}\right)  ^{2}}{2\left(  Q_{23}+Q_{34}-Q_{12}-Q_{14}\right)  }%
\end{align*}
provided the denominators are not zero.
\end{theorem}

\begin{proof}
Both $\left\{  Q_{12},Q_{23},Q_{13}\right\}  $ and $\left\{  Q_{14}%
,Q_{13},Q_{34}\right\}  $ are quad triples, so the Two quad triples theorem
shows that $Q\left(  Q_{12},Q_{23},Q_{34},Q_{14}\right)  =0$ and gives the
stated formula for $Q_{13}.$ The result for $Q_{24}$ is similar.
\end{proof}

\begin{theorem}
[Brahmagupta's identity]Suppose that $D_{12}\equiv d_{12}^{2}$, $D_{23}\equiv
d_{23}^{2}$, $D_{34}\equiv d_{34}^{2}$ and $D_{14}\equiv d_{14}^{2}$ for some
numbers $d_{12},d_{23},d_{34}$ and $d_{14}.$ Then%
\begin{align*}
&  Q\left(  D_{12},D_{23},D_{34},D_{14}\right) \\
&  =\left(  d_{12}-d_{14}+d_{23}+d_{34}\right)  \left(  d_{12}+d_{14}%
+d_{23}-d_{34}\right)  \left(  d_{14}-d_{12}+d_{23}+d_{34}\right)
\allowbreak\\
&  \times\left(  d_{12}+d_{14}-d_{23}+d_{34}\right)  \left(  d_{12}%
+d_{14}+d_{23}+d_{34}\right)  \left(  d_{12}-d_{14}-d_{23}+d_{34}\right) \\
&  \times\left(  d_{12}-d_{14}+d_{23}-d_{34}\right)  \left(  d_{23}%
-d_{14}-d_{12}+d_{34}\right)  .
\end{align*}

\end{theorem}

\begin{proof}
Make the substitutions $D_{ij}=d_{ij}^{2}$ for all $i$ and $j$ to turn the
expression%
\[
\left(  \left(  D_{12}+D_{23}+D_{34}+D_{14}\right)  ^{2}-2\left(  D_{12}%
^{2}+D_{23}^{2}+D_{34}^{2}+D_{14}^{2}\right)  \right)  ^{2}-64D_{12}%
D_{23}D_{34}D_{14}%
\]
into a difference of squares. This is then the product of the expression
\begin{align*}
&  \left(  d_{14}^{2}+d_{34}^{2}+d_{12}^{2}+d_{23}^{2}\right)  ^{2}-2\left(
d_{14}^{4}+d_{34}^{4}+d_{12}^{4}+d_{23}^{4}\right)  +8d_{14}d_{34}d_{12}%
d_{23}\\
&  =\left(  -d_{12}+d_{14}+d_{23}+d_{34}\right)  \left(  d_{12}-d_{14}%
+d_{23}+d_{34}\right) \\
&  \times\left(  d_{12}+d_{14}-d_{23}+d_{34}\right)  \left(  d_{12}%
+d_{14}+d_{23}-d_{34}\right)
\end{align*}
and the expression%
\begin{align*}
&  \left(  d_{14}^{2}+d_{34}^{2}+d_{12}^{2}+d_{23}^{2}\right)  ^{2}-2\left(
d_{14}^{4}+d_{34}^{4}+d_{12}^{4}+d_{23}^{4}\right)  -8d_{14}d_{34}d_{12}%
d_{23}\\
&  =\left(  d_{12}+d_{14}+d_{23}+d_{34}\right)  \left(  d_{12}-d_{14}%
-d_{23}+d_{34}\right) \\
&  \times\left(  d_{12}-d_{14}+d_{23}-d_{34}\right)  \left(  d_{23}%
-d_{14}-d_{12}+d_{34}\right)  .
\end{align*}

\end{proof}

By a classical theorem of Brahmagupta, the first of these two expressions
corresponds in the decimal number system to sixteen times the square of the
area of a convex cyclic quadrilateral with side lengths $d_{1},d_{2},d_{3}$
and $d_{4}.$ The second corresponds to an analogous result for a non-convex
cyclic quadrilateral with these side lengths, as discussed in \cite{Robbins}.

\subsection{Higher quad formulas}

Generalizations of the Triple quad formula and Quadruple quad formula exist
for more than four points. Writing down pleasant expressions for these seems
an interesting challenge in pure algebra.

\subsection{Isometries of the affine line}

We adopt the convention that the image of the point $A$ under the function, or
map, $\sigma$ is denoted $A\sigma,$ and that if $\sigma_{1}$ and $\sigma_{2}$
are two maps, $\sigma_{1}\sigma_{2}$ denotes the composite map given by%
\[
A\left(  \sigma_{1}\sigma_{2}\right)  =\left(  A\sigma_{1}\right)  \sigma_{2}.
\]
This allows us to write $A\left(  \sigma_{1}\sigma_{2}\right)  $ simply as
$A\sigma_{1}\sigma_{2}.$

An \textbf{isometry of the affine line} is a function $\sigma$ which inputs
and outputs points and preserves quadrance, in the sense that for any points
$A$ and $B,$%
\[
Q\left(  A,B\right)  =Q\left(  A\sigma,B\sigma\right)  .
\]

\begin{theorem}
An isometry $\sigma$ of the affine line has exactly one of the two forms:%
\[
\left[  x\right]  \sigma=\left[  x+\alpha\right]
\]
or%
\[
\left[  x\right]  \sigma=\left[  \alpha-x\right]
\]
for some number $\alpha.$
\end{theorem}

\begin{proof}
Let $O=\left[  0\right]  $ and $I=\left[  1\right]  .$ Suppose $\sigma$ is an
isometry with $O\sigma=A=\left[  \alpha\right]  $ and $I\sigma=B=\left[
\beta\right]  $. Then we must have
\[
Q\left(  A,B\right)  =\left(  \beta-\alpha\right)  ^{2}=Q\left(  O,I\right)
=1
\]
and it follows that $\beta=\alpha\pm1.$ Suppose that $\beta=\alpha+1.$ In this
case if $\left[  x\right]  \sigma=\left[  y\right]  $, then
\begin{align*}
x^{2}  & =\left(  y-\alpha\right)  ^{2}\\
\left(  x-1\right)  ^{2}  & =\left(  y-\beta\right)  ^{2}=\left(
y-\alpha-1\right)  ^{2}.
\end{align*}
It follows that $y=x+\alpha.$ Now suppose that $\beta=\alpha-1.$ In this case
if $\left[  x\right]  \sigma=\left[  y\right]  $, then
\begin{align*}
x^{2}  & =\left(  y-\alpha\right)  ^{2}\\
\left(  x-1\right)  ^{2}  & =\left(  y-\beta\right)  ^{2}=\left(
y-\alpha+1\right)  ^{2}.
\end{align*}
It follows that $y=\alpha-x.$ Each of these maps is easily checked to be an isometry.
\end{proof}

Note in particular that an isometry is invertible.

\section{The projective line}

\subsection{Projective quadrance}

A \textbf{projective point} is a proportion $a\equiv\left[  x:y\right]  $
where $x$ and $y$ are not both zero, and where $\left[  x_{1}:y_{1}\right]
=\left[  x_{2}:y_{2}\right]  $ precisely when
\[
x_{1}y_{2}-x_{2}y_{1}=0.
\]
This is equivalent to the condition $\left[  x:y\right]  =\left[  \lambda
x:\lambda y\right]  $ for any non-zero $\lambda.$ A projective point will
often be called a \textbf{p-point}.

A \textbf{form} is a proportion $F\equiv\left(  d:e:f\right)  $ where $d,e$
and $f$ are not all zero. Again we agree that
\[
\left(  d:e:f\right)  =\left(  \lambda d:\lambda e:\lambda f\right)
\]
for any non-zero $\lambda.$ The form $F\equiv\left(  d:e:f\right)  $ is
\textbf{degenerate }precisely when the \textbf{discriminant}
\[
df-e^{2}%
\]
is zero. Everything in this section depends on the choice of a non-degenerate
form $F\equiv\left(  d:e:f\right)  $, which we consider arbitrary but fixed.

A p-point $a\equiv\left[  x:y\right]  $ is \textbf{null} precisely when%
\[
dx^{2}+2exy+fy^{2}=0.
\]
The p-points $a_{1}\equiv\left[  x_{1}:y_{1}\right]  $ and $a_{2}\equiv\left[
x_{2}:y_{2}\right]  $ are \textbf{perpendicular} precisely when
\[
dx_{1}x_{2}+ex_{1}y_{2}+ex_{2}y_{1}+fy_{1}y_{2}=0.
\]

The \textbf{projective quadrance}, or more simply the \textbf{p-quadrance},
between non-null p-points $a_{1}\equiv\left[  x_{1}:y_{1}\right]  $ and
$a_{2}\equiv\left[  x_{2}:y_{2}\right]  $ is the number%
\[
q=q\left(  a_{1},a_{2}\right)  \equiv\frac{\left(  df-e^{2}\right)  \left(
x_{1}y_{2}-x_{2}y_{1}\right)  ^{2}}{\left(  dx_{1}^{2}+2ex_{1}y_{1}+fy_{1}%
^{2}\right)  \left(  dx_{2}^{2}+2ex_{2}y_{2}+fy_{2}^{2}\right)  }.
\]
The presence of the discriminant $df-e^{2}$ ensures that the formulas that
follow are independent of the form $F.$ Note that
\[
q\left(  a_{1},a_{2}\right)  =q\left(  a_{2},a_{1}\right)
\]
and that $q\left(  a_{1},a_{2}\right)  $ is zero precisely when $a_{1}=a_{2}.$

\begin{theorem}
[Perpendicular p-points]Two non-null p-points $a_{1}$ and $a_{2}$ are
perpendicular precisely when $q\left(  a_{1},a_{2}\right)  =1.$
\end{theorem}

\begin{proof}
The \textbf{generalized Fibonacci's identity}%
\begin{align*}
& \left(  df-e^{2}\right)  \left(  x_{1}y_{2}-x_{2}y_{1}\right)  ^{2}+\left(
dx_{1}x_{2}+ex_{1}y_{2}+ex_{2}y_{1}+fy_{1}y_{2}\right)  ^{2}\\
& =\left(  dx_{1}^{2}+2ex_{1}y_{1}+fy_{1}^{2}\right)  \left(  dx_{2}%
^{2}+2ex_{2}y_{2}+fy_{2}^{2}\right)
\end{align*}
shows that if $a_{1}\equiv\left[  x_{1}:y_{1}\right]  $ and $a_{2}%
\equiv\left[  x_{2}:y_{2}\right]  $, then%
\[
1-q\left(  a_{1},a_{2}\right)  =\frac{\left(  dx_{1}x_{2}+ex_{1}y_{2}%
+ex_{2}y_{1}+fy_{1}y_{2}\right)  ^{2}}{\left(  dx_{1}^{2}+2ex_{1}y_{1}%
+fy_{1}^{2}\right)  \left(  dx_{2}^{2}+2ex_{2}y_{2}+fy_{2}^{2}\right)  }.
\]
This is zero precisely when $a_{1}$ and $a_{2}$ are perpendicular.
\end{proof}

The next result is fundamental for metrical projective geometry. It is the
analog of the Triple quad formula, and differs from it only by a cubic term.
In planar rational trigonometry this result is called the \textit{Triple
spread formula}.

\begin{theorem}
[Projective triple quad formula]For p-points $a_{1},a_{2}$ and $a_{3},$ define
the p-quadrances $q_{1}\equiv q\left(  a_{2},a_{3}\right)  $, $q_{2}\equiv
q\left(  a_{1},a_{3}\right)  $ and $q_{3}\equiv q\left(  a_{1},a_{2}\right)
$. Then
\[
\left(  q_{1}+q_{2}+q_{3}\right)  ^{2}=2\left(  q_{1}^{2}+q_{2}^{2}+q_{3}%
^{2}\right)  +4q_{1}q_{2}q_{3}.
\]

\end{theorem}

\begin{proof}
If $a_{1}\equiv\left[  x_{1}:y_{1}\right]  $, $a_{2}\equiv\left[  x_{2}%
:y_{2}\right]  ,$ and $a_{3}\equiv\left[  x_{3}:y_{3}\right]  $ then%
\begin{align*}
q_{1}  & =\frac{\left(  df-e^{2}\right)  \left(  x_{2}y_{3}-x_{3}y_{2}\right)
^{2}}{\left(  dx_{2}^{2}+2ex_{2}y_{2}+fy_{2}^{2}\right)  \left(  dx_{3}%
^{2}+2ex_{3}y_{3}+fy_{3}^{2}\right)  }\\
q_{2}  & =\frac{\left(  df-e^{2}\right)  \left(  x_{1}y_{3}-x_{3}y_{1}\right)
^{2}}{\left(  dx_{1}^{2}+2ex_{1}y_{1}+fy_{1}^{2}\right)  \left(  dx_{3}%
^{2}+2ex_{3}y_{3}+fy_{3}^{2}\right)  }\\
q_{3}  & =\frac{\left(  df-e^{2}\right)  \left(  x_{1}y_{2}-x_{2}y_{1}\right)
^{2}}{\left(  dx_{1}^{2}+2ex_{1}y_{1}+fy_{1}^{2}\right)  \left(  dx_{2}%
^{2}+2ex_{2}y_{2}+fy_{2}^{2}\right)  }.
\end{align*}

The following is an algebraic identity:%
\begin{align*}
& \left(  x_{2}y_{3}-x_{3}y_{2}\right)  ^{2}\left(  dx_{1}^{2}+2ex_{1}%
y_{1}+fy_{1}^{2}\right)  +\left(  x_{1}y_{3}-x_{3}y_{1}\right)  ^{2}\left(
dx_{2}^{2}+2ex_{2}y_{2}+fy_{2}^{2}\right) \\
& -\left(  x_{1}y_{2}-x_{2}y_{1}\right)  ^{2}\left(  dx_{3}^{2}+2ex_{3}%
y_{3}+fy_{3}^{2}\right) \\
& =\allowbreak2\left(  x_{2}y_{3}-x_{3}y_{2}\right)  \left(  x_{1}y_{3}%
-x_{3}y_{1}\right)  \left(  dx_{1}x_{2}+ex_{1}y_{2}+ex_{2}y_{1}+fy_{1}%
y_{2}\right)
\end{align*}
Square both sides and rearrange using the Perpendicular p-points theorem, to
get%
\[
\left(  q_{1}+q_{2}-q_{3}\right)  ^{2}=4q_{1}q_{2}\left(  1-q_{3}\right)  .
\]
This can be rewritten as the symmetrical equation
\[
\left(  q_{1}+q_{2}+q_{3}\right)  ^{2}-2\left(  q_{1}^{2}+q_{2}^{2}+q_{3}%
^{2}\right)  =4q_{1}q_{2}q_{3}.
\]

\end{proof}

Motivated by the Projective triple quad formula, define the \textbf{Triple
spread function }$S\left(  a,b,c\right)  $ for numbers $a,b$ and $c$ by%
\[
S\left(  a,b,c\right)  \equiv\left(  a+b+c\right)  ^{2}-2\left(  a^{2}%
+b^{2}+c^{2}\right)  -4abc.
\]
The reason for the terminology is that in two dimensional geometry the
projective quadrance between p-points becomes the \textit{spread} between
lines. Note that $S\left(  a,b,c\right)  $ is a symmetric function of $a,b$
and $c,$ and that%
\begin{align}
S\left(  a,b,c\right)   & =A\left(  a,b,c\right)  -4abc\nonumber\\
& =2\left(  ac+bc+ab\right)  -\left(  a^{2}+b^{2}+c^{2}\right)
-4abc\nonumber\\
& =4\left(  ab+bc+ca\right)  -\left(  a+b+c\right)  ^{2}-4abc\nonumber\\
& =4\left(  1-a\right)  \left(  1-b\right)  \left(  1-c\right)  -\left(
a+b+c-2\right)  ^{2}\nonumber\\
& =4\left(  1-a\right)  bc-\left(  a-b-c\right)  ^{2}\nonumber\\
& =-%
\begin{vmatrix}
0 & a & b & 1\\
a & 0 & c & 1\\
b & c & 0 & 1\\
1 & 1 & 1 & 2
\end{vmatrix}
\nonumber\\
& =4bc\left(  1-b\right)  \left(  1-c\right)  -\left(  a-b-c+2bc\right)
^{2}.\label{TripleSpread}%
\end{align}
\qquad The last of these equations is useful to solve $S\left(  a,b,c\right)
=0$ for $a$ if $b$ and $c$ are known.

A set $\left\{  a,b,c\right\}  $ is a \textbf{spread triple} precisely when
$S\left(  a,b,c\right)  =0.$

\subsection{Projective quadruple quad formula}

The next theorems extend the Projective triple quad formula to four p-points.

\begin{theorem}
[Two spread triples]Suppose that $\left\{  a,b,x\right\}  $ and $\left\{
c,d,x\right\}  $ are both spread triples. Then%
\begin{align*}
&  \left(
\begin{array}
[c]{c}%
\left(  a+b+c+d\right)  ^{2}-2\left(  a^{2}+b^{2}+c^{2}+d^{2}\right) \\
-4\left(  abc+abd+acd+bcd\right)  +8abcd
\end{array}
\right)  ^{2}\\
&  =64abcd\left(  1-a\right)  \left(  1-b\right)  \left(  1-c\right)  \left(
1-d\right)  .
\end{align*}
Furthermore if $a+b-2ab\neq c+d-2cd$ then
\[
x=\frac{\left(  a-b\right)  ^{2}-\left(  c-d\right)  ^{2}}{2\left(
a+b-c-d-2ab+2cd\right)  }.
\]

\end{theorem}

\begin{proof}
Suppose that $\left\{  a,b,x\right\}  $ and $\left\{  c,d,x\right\}  $ are
both spread triples. Then (\ref{TripleSpread}) gives
\begin{align}
\left(  x-a-b+2ab\right)  ^{2} &  =4ab\left(  1-a\right)  \left(  1-b\right)
\label{TripleSp1}\\
\left(  x-c-d+2cd\right)  ^{2} &  =4cd\left(  1-c\right)  \left(  1-d\right)
.\label{TripleSp2}%
\end{align}
Take the difference between these two equations. If $a+b-2ab\neq c+d-2cd$ then
you may solve for $x$ to get%
\[
x=\frac{\left(  a-b\right)  ^{2}-\left(  c-d\right)  ^{2}}{2\left(
a+b-c-d-2ab+2cd\right)  }.
\]
Substitute back into (\ref{TripleSp1}) or (\ref{TripleSp2}) to get%
\begin{align}
& \left(  \left(  a-b\right)  ^{2}-\left(  c-d\right)  ^{2}-2\left(
a+b-c-d-2ab+2cd\right)  \left(  a+b-2ab\right)  \right)  ^{2}\nonumber\\
& =16ab\left(  1-a\right)  \left(  1-b\right)  \left(  a+b-c-d-2ab+2cd\right)
^{2}.\label{QuadQuad}%
\end{align}
$\allowbreak$This condition can be rewritten more symmetrically as
\begin{align*}
&  \left(
\begin{array}
[c]{c}%
\left(  a+b+c+d\right)  ^{2}-2\left(  a^{2}+b^{2}+c^{2}+d^{2}\right) \\
-4\left(  abc+abd+acd+bcd\right)  +8abcd
\end{array}
\right)  ^{2}\\
&  =64abcd\left(  1-a\right)  \left(  1-b\right)  \left(  1-c\right)  \left(
1-d\right)  .
\end{align*}
If
\[
a+b-2ab=c+d-2cd
\]
then by (\ref{TripleSp1}) also
\[
4ab\left(  1-a\right)  \left(  1-b\right)  =4cd\left(  1-c\right)  \left(
1-d\right)  .
\]
$\allowbreak$ Hence the identity%
\[
\left(  a+b-2ab\right)  ^{2}-4ab\left(  1-a\right)  \left(  1-b\right)
=\allowbreak\left(  a-b\right)  ^{2}%
\]
implies that $\left(  a-b\right)  ^{2}=\left(  c-d\right)  ^{2},$ and then
(\ref{QuadQuad}) is automatic.
\end{proof}

The \textbf{Quadruple spread function }$R\left(  a,b,c,d\right)  $ is defined
for numbers $a,b,c$ and $d$ by%
\begin{align*}
R\left(  a,b,c,d\right)   &  \equiv\left(
\begin{array}
[c]{c}%
\left(  a+b+c+d\right)  ^{2}-2\left(  a^{2}+b^{2}+c^{2}+d^{2}\right) \\
-4\left(  abc+abd+acd+bcd\right)  +8abcd
\end{array}
\right)  ^{2}\\
&  -64abcd\left(  1-a\right)  \left(  1-b\right)  \left(  1-c\right)  \left(
1-d\right)  .
\end{align*}
This is a symmetric function of $a,b,c$ and $d$.

\begin{theorem}
[Projective quadruple quad formula]Suppose $a_{1},a_{2},a_{3}$ and $a_{4}$ are
non-null p-points with p-quadrances $q_{ij}\equiv q\left(  a_{i},a_{j}\right)
$ for all $i,j=1,2,3$ and $4$. Then%
\[
R\left(  q_{12},q_{23},q_{34},q_{14}\right)  =0.
\]
Furthermore%
\begin{align*}
q_{13} &  =\frac{\left(  q_{12}-q_{23}\right)  ^{2}-\left(  q_{34}%
-q_{14}\right)  ^{2}}{2\left(  q_{12}+q_{23}-q_{34}-q_{14}-2q_{12}%
q_{23}+2q_{34}q_{14}\right)  }\\
q_{24} &  =\frac{\left(  q_{23}-q_{34}\right)  ^{2}-\left(  q_{12}%
-q_{14}\right)  ^{2}}{2\left(  q_{23}+q_{34}-q_{12}-q_{14}-2q_{23}%
q_{34}+2q_{12}q_{14}\right)  }%
\end{align*}
provided the denominators are non-zero.
\end{theorem}

\begin{proof}
Since $\left\{  q_{12},q_{23},q_{13}\right\}  $ and $\left\{  q_{14}%
,q_{34},q_{13}\right\}  $ are both spread triples, and $\left\{  q_{23}%
,q_{34},q_{24}\right\}  $ and $\left\{  q_{12},q_{14},q_{24}\right\}  $ are
also both spread triples, the formulas follow from the Two spread triples theorem.
\end{proof}

\subsection{An example}

Suppose $a_{1}\equiv\left[  1:0\right]  $, $a_{2}\equiv\left[  2:3\right]  $,
$a_{3}\equiv\left[  4:-1\right]  $ and $a_{4}\equiv\left[  3:5\right]  $ over
the rational number field, and that we choose the Euclidean form
$F\equiv\left(  1:0:1\right)  $, corresponding to $x^{2}+y^{2}.$ Then the
projective quadrance is given by
\[
q\left(  \left[  x_{1}:y_{1}\right]  ,\left[  x_{2}:y_{2}\right]  \right)
=\frac{\left(  x_{1}y_{2}-x_{2}y_{1}\right)  ^{2}}{\left(  x_{1}^{2}+y_{1}%
^{2}\right)  \left(  x_{2}^{2}+y_{2}^{2}\right)  }.
\]

You may then calculate that with $q_{ij}\equiv q\left(  a_{i},a_{j}\right)
$,
\[%
\begin{tabular}
[c]{lllllll}%
$q_{12}=9/13$ &  & $q_{23}=196/221$ &  & $q_{34}=529/578$ &  & $q_{14}=25/34$%
\end{tabular}
\]
and that%

\[
R\left(  \frac{9}{13},\frac{196}{221},\allowbreak\frac{529}{578},\frac{25}%
{34}\right)  =0.
\]
Furthermore%

\begin{align*}
\frac{\left(  q_{12}-q_{23}\right)  ^{2}-\left(  q_{34}-q_{14}\right)  ^{2}%
}{2\left(  q_{12}+q_{23}-q_{34}-q_{14}-2q_{12}q_{23}+2q_{34}q_{14}\right)  }
& =\allowbreak\frac{1}{17}=q_{13}\\
\frac{\left(  q_{23}-q_{34}\right)  ^{2}-\left(  q_{12}-q_{14}\right)  ^{2}%
}{2\left(  q_{23}+q_{34}-q_{12}-q_{14}-2q_{23}q_{34}+2q_{12}q_{14}\right)  }
& =\allowbreak\frac{1}{442}=q_{24}.
\end{align*}

\subsection{Higher quadruple spread formulas}

The same remarks made earlier about higher quad formulas apply to higher
spread formulas.

\section{The Spread polynomials}

The Projective triple quad formula, or Triple spread formula%
\[
\left(  a+b+c\right)  ^{2}=2\left(  a^{2}+b^{2}+c^{2}\right)  +4abc
\]
yields an important sequence of polynomials with integer coefficients, called
\textit{spread polynomials}. These are universal analogs of the Chebyshev
polynomials of the first kind, and make sense over any field, not of
characteristic two. It is interesting that these appear already in one
dimensional geometry.

To motivate the spread polynomials, note that in the Triple spread formula if
$a\equiv b\equiv s$ then $c$ is either $0$ or $4s\left(  1-s\right)  $. If
$a\equiv4s\left(  1-s\right)  $ and $b\equiv s$ then $c$ is either $s$ or
$s\left(  3-4s\right)  ^{2}.$ Continuing, there is a sequence of polynomials
$S_{n}\left(  s\right)  $ for $n=0,1,2,\cdots$ with the property that
$S_{n-1}\left(  s\right)  ,s$ and $S_{n}\left(  s\right)  $ always satisfy the
Triple spread formula.

As in \cite{Wild}, the \textbf{spread polynomial }$S_{n}\left(  s\right)  $ is
defined recursively by $S_{0}\left(  s\right)  \equiv0$ and $S_{1}\left(
s\right)  \equiv s,$ together with the rule%
\[
S_{n}\left(  s\right)  \equiv2\left(  1-2s\right)  S_{n-1}\left(  s\right)
-S_{n-2}\left(  s\right)  +2s.
\]

The coefficient of $s^{n}$ in $S_{n}\left(  s\right)  $ is a power of four, so
in any field not of characteristic two the degree of the polynomial
$S_{n}\left(  s\right)  $ is $n$. Over the decimal number field
\[
S_{n}\left(  s\right)  =\frac{1-T_{n}\left(  1-2s\right)  }{2}%
\]
where $T_{n}$ is the $n$-th Chebyshev polynomial of the first kind. The affine
transformation which maps the square $\left[  -1,1\right]  \times\left[
-1,1\right]  $ to the square $\left[  0,1\right]  \times\left[  0,1\right]  $
takes the graph of the $n$-th Chebyshev polynomial to the graph of the $n$-th
spread polynomial.

The first few spread polynomials are%
\begin{align*}
S_{0}\left(  s\right)   & =0\\
S_{1}\left(  s\right)   & =s\\
S_{2}\left(  s\right)   & =4s-4s^{2}=4s\left(  1-s\right) \\
S_{3}\left(  s\right)   & =9s-24s^{2}+16s^{3}=\allowbreak s\left(
4s-3\right)  ^{2}\\
S_{4}\left(  s\right)   & =16s-80s^{2}+128s^{3}-64s^{4}=16s\left(  1-s\right)
\left(  2s-1\right)  ^{2}\\
S_{5}\left(  s\right)   & =25s-200s^{2}+560s^{3}-640s^{4}+256s^{5}=\allowbreak
s\left(  16s^{2}-20s+5\right)  ^{2}.
\end{align*}
Note that $S_{2}\left(  s\right)  $ is the logistic map. As shown in
\cite{Wild}, $S_{n}\circ S_{m}=S_{nm}$ for $n,m\geq1,$ and the spread
polynomials have interesting orthogonality properties over finite fields.

S. Goh \cite{Goh} observed that there is a sequence of `spread-cyclotomic'
polynomials $\phi_{k}\left(  s\right)  $ of degree $\phi\left(  k\right)  $
with integer coefficients such that for any $n=1,2,3,\cdots$
\[
S_{n}\left(  s\right)  =\prod_{k|n}\phi_{k}\left(  s\right)  .
\]
This factorization shows that arithmetically the spread polynomials have quite
different properties than the closely related Chebyshev polynomials.

\section{Chromogeometry}

There are three forms which are particularly important and useful, and which
interact in a surprising way. This becomes particularly clear in two
dimensional geometry, but there are already hints here in the one dimensional situation.

The \textbf{blue form} is $F_{b}\equiv\left(  1:0:1\right)  $, the \textbf{red
form} is $F_{r}\equiv\left(  1:0:-1\right)  $ and the \textbf{green form} is
$F_{g}\equiv\left(  0:1:0\right)  $. Two p-points $a_{1}\equiv\left[
x_{1}:y_{1}\right]  $ and $a_{2}\equiv\left[  x_{2}:y_{2}\right]  $ are then
\textbf{blue, red and green perpendicular} precisely when%
\begin{align*}
x_{1}x_{2}+y_{1}y_{2}  & =0\\
x_{1}x_{2}-y_{1}y_{2}  & =0\\
x_{1}y_{2}+x_{2}y_{1}  & =0
\end{align*}
respectively. Recall that $a_{1}=a_{2}$ precisely when%
\[
x_{1}y_{2}-x_{2}y_{1}=0
\]
and these four equations exhaust the bilinear equations in the coordinates of
$a_{1}$ and $a_{2}$ which employ coefficients of $\pm1.$

The \textbf{blue, red and green perpendiculars} of the p-point $a\equiv\left[
x:y\right]  $ are respectively the p-points
\[%
\begin{tabular}
[c]{lllll}%
$a^{b}\equiv\left[  -y:x\right]  $ &  & $a^{r}\equiv\left[  y:x\right]  $ &  &
$a^{g}\equiv\left[  x:-y\right]  .$%
\end{tabular}
\]

We let $q^{b},q^{r}$ and $q^{g}$ denote the blue, red and green p-quadrances
associated to the blue, red and green forms respectively. Then%
\begin{align*}
q^{b}\left(  \left[  x_{1}:y_{1}\right]  ,\left[  x_{2}:y_{2}\right]  \right)
& =\frac{\left(  x_{1}y_{2}-x_{2}y_{1}\right)  ^{2}}{\left(  x_{1}^{2}%
+y_{1}^{2}\right)  \left(  x_{2}^{2}+y_{2}^{2}\right)  }\\
q^{r}\left(  \left[  x_{1}:y_{1}\right]  ,\left[  x_{2}:y_{2}\right]  \right)
& =-\frac{\left(  x_{1}y_{2}-x_{2}y_{1}\right)  ^{2}}{\left(  x_{1}^{2}%
-y_{1}^{2}\right)  \left(  x_{2}^{2}-y_{2}^{2}\right)  }\\
q^{g}\left(  \left[  x_{1}:y_{1}\right]  ,\left[  x_{2}:y_{2}\right]  \right)
& =-\frac{\left(  x_{1}y_{2}-x_{2}y_{1}\right)  ^{2}}{4x_{1}y_{1}x_{2}y_{2}}.
\end{align*}

\begin{theorem}
For any p-points $a_{1}$ and $a_{2}$
\[
\frac{1}{q^{b}\left(  a_{1},a_{2}\right)  }+\frac{1}{q^{r}\left(  a_{1}%
,a_{2}\right)  }+\frac{1}{q^{g}\left(  a_{1},a_{2}\right)  }=2.
\]

\end{theorem}

\begin{proof}
This follows from the previous formulas and the identity%
\[
\left(  x_{1}^{2}+y_{1}^{2}\right)  \left(  x_{2}^{2}+y_{2}^{2}\right)
-\left(  x_{1}^{2}-y_{1}^{2}\right)  \left(  x_{2}^{2}-y_{2}^{2}\right)
-4x_{1}y_{1}x_{2}y_{2}=\allowbreak2\left(  x_{1}y_{2}-x_{2}y_{1}\right)  ^{2}.
\]

\end{proof}

\begin{theorem}
For any p-point $a,$
\[
q^{b}\left(  a^{r},a^{g}\right)  =q^{r}\left(  a^{g},a^{b}\right)
=q^{g}\left(  a^{b},a^{r}\right)  =1.
\]

\end{theorem}

\begin{proof}
For $a\equiv\left[  x:y\right]  $
\begin{align*}
q^{b}\left(  a^{r},a^{g}\right)   & =q^{b}\left(  \left[  y:x\right]  ,\left[
x:-y\right]  \right)  =1\\
q^{r}\left(  a^{g},a^{b}\right)   & =q^{r}\left(  \left[  x:-y\right]
,\left[  -y:x\right]  \right)  =1\\
q^{g}\left(  a^{b},a^{r}\right)   & =q^{g}\left(  \left[  -y:x\right]
,\left[  y:x\right]  \right)  =1.
\end{align*}

\end{proof}

\begin{theorem}
For any two p-points $a_{1}=\left[  x_{1}:y_{1}\right]  $ and $a_{2}=\left[
x_{2}:y_{2}\right]  $ and any colour $c$ (either $b,r$ or $g$)%
\[
q^{c}\left(  a_{1},a_{2}\right)  =q^{c}\left(  a_{1}^{b},a_{2}^{b}\right)
=q^{c}\left(  a_{1}^{r},a_{2}^{r}\right)  =q^{c}\left(  a_{1}^{g},a_{2}%
^{g}\right)  .
\]

\end{theorem}

\begin{proof}
A computation using the above formulas.
\end{proof}

\begin{theorem}
For any p-points $a_{1}$ and $a_{2},$%
\begin{align*}
q^{b}\left(  a_{1}^{r},a_{2}^{g}\right)   & =q^{b}\left(  a_{1}^{g},a_{2}%
^{r}\right) \\
q^{r}\left(  a_{1}^{g},a_{2}^{b}\right)   & =q^{r}\left(  a_{1}^{b},a_{2}%
^{g}\right) \\
q^{g}\left(  a_{1}^{b},a_{2}^{r}\right)   & =q^{g}\left(  a_{1}^{r},a_{2}%
^{b}\right)  .
\end{align*}

\end{theorem}

\begin{proof}
A computation.
\end{proof}

\section{Projective isometries: blue, red and green}

Let $q$ denote one of the blue, red or green p-quadrances. A
projective\textbf{\ isometry} is a map $\sigma$ that inputs and outputs
projective points and satisfies%
\[
q\left(  a,b\right)  =q\left(  a\sigma,b\sigma\right)
\]
for any projective points $a$ and $b.$ Note that the action of the map
$\sigma$ on the p-point $a$ is denote $a\sigma.$

One way of describing a map on projective points is to use a
\textbf{projective matrix}%
\[%
\begin{bmatrix}
a & b\\
c & d
\end{bmatrix}
\]
where by definition%
\[%
\begin{bmatrix}
a & b\\
c & d
\end{bmatrix}
=%
\begin{bmatrix}
\lambda a & \lambda b\\
\lambda c & \lambda d
\end{bmatrix}
\]
for any non-zero number $\lambda.$ Such a matrix defines the map
\[
\left[  x:y\right]  \sigma=\left[  x:y\right]
\begin{bmatrix}
a & b\\
c & d
\end{bmatrix}
=\left[  ax+cy:bx+dy\right]
\]
and we write%
\[
\sigma=%
\begin{bmatrix}
a & b\\
c & d
\end{bmatrix}
.
\]

\subsection{Blue isometries}

\begin{theorem}
An isometry of the blue p-quadrance is either%
\[%
\begin{tabular}
[c]{lllll}%
$\sigma_{\left[  a:b\right]  }^{b}=%
\begin{bmatrix}
a & b\\
b & -a
\end{bmatrix}
$ &  & \textrm{or} &  & $\rho_{\left[  a:b\right]  }^{b}=%
\begin{bmatrix}
a & b\\
-b & a
\end{bmatrix}
$%
\end{tabular}
\]
for some blue non-null projective point $\left[  a:b\right]  $.
\end{theorem}

\begin{proof}
Suppose that $\sigma$ sends $i_{1}\equiv\left[  1:0\right]  $ to $\left[
a:b\right]  $ and $i_{2}\equiv\left[  0:1\right]  $ to $\left[  c:d\right]  $.
Then since $q_{b}\left(  i_{1},i_{2}\right)  =1,$ we must have
\[
ac+bd=0.
\]
So $\left[  c:d\right]  =\left[  b:-a\right]  .$ Now given an arbitrary
projective point $u\equiv\left[  x:y\right]  $ with $u\sigma\equiv
v\equiv\left[  w:z\right]  $,%
\begin{align*}
q_{b}\left(  u,i_{1}\right)   & =1-\frac{x^{2}}{x^{2}+y^{2}}=\frac{y^{2}%
}{x^{2}+y^{2}}\\
& =q\left(  v,\left[  a:b\right]  \right)  =\frac{\left(  az-bw\right)  ^{2}%
}{\left(  a^{2}+b^{2}\right)  \left(  w^{2}+z^{2}\right)  }%
\end{align*}
and%
\begin{align*}
q_{b}\left(  u,i_{2}\right)   & =1-\frac{y^{2}}{x^{2}+y^{2}}=\frac{x^{2}%
}{x^{2}+y^{2}}\\
& =q\left(  v,\left[  b:-a\right]  \right)  =\frac{\left(  bz+aw\right)  ^{2}%
}{\left(  a^{2}+b^{2}\right)  \left(  w^{2}+z^{2}\right)  }.
\end{align*}
This gives the following quadratic equations for $w$ and $z:$
\begin{align*}
\frac{y^{2}}{x^{2}+y^{2}}  & =\frac{\left(  az-bw\right)  ^{2}}{\left(
a^{2}+b^{2}\right)  \left(  w^{2}+z^{2}\right)  }\\
\frac{x^{2}}{x^{2}+y^{2}}  & =\frac{\left(  bz+aw\right)  ^{2}}{\left(
a^{2}+b^{2}\right)  \left(  w^{2}+z^{2}\right)  }.
\end{align*}
The two possible solutions for $\left[  w:z\right]  $ are%
\[%
\begin{tabular}
[c]{lllll}%
$\left[  ax+by:bx-ay\right]  $ &  & \textrm{and} &  & $\left[
ax-by:bx+ay\right]  $%
\end{tabular}
\]
which correspond to $u\sigma_{\left[  a:b\right]  }^{b}$ and $u\rho_{\left[
a:b\right]  }^{b}$ respectively.
\end{proof}

\begin{theorem}
For any blue non-null projective points $\left[  a:b\right]  $ and $\left[
c:d\right]  $
\end{theorem}%

\[%
\begin{tabular}
[c]{lll}%
$\sigma_{\left[  a:b\right]  }^{b}\sigma_{\left[  c:d\right]  }^{b}%
=\rho_{\left[  ac+bd:ad-bc\right]  }^{b}$ &  & $\rho_{\left[  a:b\right]
}^{b}\rho_{\left[  c:d\right]  }^{b}=\rho_{\left[  ac-bd:ad+bc\right]  }^{b}%
$\\
&  & \\
$\rho_{\left[  a:b\right]  }^{b}\sigma_{\left[  c:d\right]  }^{b}%
=\sigma_{\left[  ac+bd:ad-bc\right]  }^{b}$ &  & $\sigma_{\left[  a:b\right]
}^{b}\rho_{\left[  c:d\right]  }^{b}=\sigma_{\left[  ac-bd:ad+bc\right]  }%
^{b}.$%
\end{tabular}
\]

\begin{proof}
This is a straightforward verification. Note that the Fibonacci identities%
\[
\left(  ac+bd\right)  ^{2}+\left(  ad-bc\right)  ^{2}=\left(  a^{2}%
+b^{2}\right)  \left(  c^{2}+d^{2}\right)  =\left(  ac-bd\right)  ^{2}+\left(
ad+bc\right)  ^{2}%
\]
show that the resultant isometries are also associated to non-null points.
\end{proof}

Define the group $G^{b}$ of all isometries of the blue projective line, and
distinguish the subgroup $G_{e}^{b}$ of \textbf{blue rotations} $\rho_{\left[
a:b\right]  }^{b}.$ These latter are naturally in bijection with the non-null
projective points.

The coset $G_{o}^{b}$ consists of \textbf{blue reflections} $\sigma_{\left[
a:b\right]  }^{b}$ which are also naturally in bijection with the non-null
points. The non-commutative group $G^{b}$ naturally acts on the blue non-null
projective points, as well as on the blue null projective points. Both actions
are transitive.

The group structure of the blue rotations can be transferred to give a
\textbf{multiplication} of blue non-null projective points, by the rule%
\[
\left[  a:b\right]  \times_{b}\left[  c:d\right]  =\left[  ac-bd:ad+bc\right]
.
\]
This multiplication is associative, commutative, has identity $\left[
1:0\right]  ,$ and the inverse of $\left[  a:b\right]  $ is $\left[
a:-b\right]  $. Of course this rule anticipates multiplication of unit complex
numbers, but there are subtle differences. For any p-point $\left[
a:b\right]  $ with $a^{2}+b^{2}=1$ and $a\neq-1$ you can compute that%
\begin{align*}
\left[  a+1:b\right]  \times_{b}\left[  a+1:b\right]   & =\left[  \left(
a+1\right)  ^{2}-b^{2}:2\left(  a+1\right)  b\right] \\
& =\left[  2\left(  a^{2}+a\right)  :2\left(  a+1\right)  b\right]  =\left[
a:b\right]
\end{align*}
so that such a p-point has a distinguished algebraic `square root' without
requiring a field extension.

\subsection{Red Isometries}

\begin{theorem}
An isometry of the red p-quadrance is either%
\[%
\begin{tabular}
[c]{lllll}%
$\sigma_{\left[  a:b\right]  }^{r}=%
\begin{bmatrix}
a & b\\
-b & -a
\end{bmatrix}
$ &  & \textrm{or} &  & $\rho_{\left[  a:b\right]  }^{r}=%
\begin{bmatrix}
a & b\\
b & a
\end{bmatrix}
$%
\end{tabular}
\]
for some red non-null projective point $\left[  a:b\right]  $.
\end{theorem}

\begin{proof}
Suppose that $\sigma$ sends $i_{1}\equiv\left[  1:0\right]  $ to $\left[
a:b\right]  $ and $i_{2}\equiv\left[  0:1\right]  $ to $\left[  c:d\right]  $.
Now given an arbitrary p-point $u\equiv\left[  x:y\right]  $ with
$u\sigma\equiv\left[  w:z\right]  $, a similar analysis to the previous proof
gives the two quadratic equations
\begin{align*}
\frac{y^{2}}{x^{2}-y^{2}}  & =\frac{\left(  az-bw\right)  ^{2}}{\left(
a^{2}-b^{2}\right)  \left(  w^{2}-z^{2}\right)  }\\
\frac{x^{2}}{x^{2}-y^{2}}  & =\frac{\left(  bz-aw\right)  ^{2}}{\left(
a^{2}-b^{2}\right)  \left(  w^{2}-z^{2}\right)  }.
\end{align*}
The two possible solutions for $\left[  w:z\right]  $ are%
\[%
\begin{tabular}
[c]{lllll}%
$\left[  ax-by:bx-ay\right]  $ &  & \textrm{and} &  & $\left[
ax+by:bx+ay\right]  $%
\end{tabular}
\]
which correspond to $u\sigma_{\left[  a:b\right]  }$ and $u\rho_{\left[
a:b\right]  }$ respectively.
\end{proof}

\begin{theorem}
For any red non-null projective points $\left[  a:b\right]  $ and $\left[
c:d\right]  $%
\[%
\begin{tabular}
[c]{lll}%
$\sigma_{\left[  a:b\right]  }^{r}\sigma_{\left[  c:d\right]  }^{r}%
=\rho_{\left[  ac-bd:ad-bc\right]  }^{r}$ &  & $\rho_{\left[  a:b\right]
}^{r}\rho_{\left[  c:d\right]  }^{r}=\rho_{\left[  ac+bd:ad+bc\right]  }^{r}%
$\\
&  & \\
$\rho_{\left[  a:b\right]  }^{r}\sigma_{\left[  c:d\right]  }^{r}%
=\sigma_{\left[  ac-bd:ad-bc\right]  }^{r}$ &  & $\rho_{\left[  c:d\right]
}^{r}\sigma_{\left[  a:b\right]  }^{r}=\sigma_{\left[  ac+bd:ad+bc\right]
}^{r}.$%
\end{tabular}
\]

\end{theorem}

\begin{proof}
A straightforward verification. Note that the hyperbolic versions of
Fibonacci's identities:
\[
\left(  ac-bd\right)  ^{2}-\left(  ad-bc\right)  ^{2}=\left(  a^{2}%
-b^{2}\right)  \left(  c^{2}-d^{2}\right)  =\left(  ac+bd\right)  ^{2}-\left(
ad+bc\right)  ^{2}%
\]
show that the resultant isometries are also associated to red non-null points.
\end{proof}

We now define the group $G^{r}$ of all isometries of the red projective line,
and distinguish the subgroup $G_{e}^{r}=G_{e}^{r}$ of red rotations
$\rho_{\left[  a:b\right]  }^{r}.$ These latter are naturally in bijection
with the non-null projective points.

The coset $G_{o}^{r}=G_{o}^{r}$ consists of red reflections $\sigma_{\left[
a:b\right]  }^{r}$ which are also naturally in bijection with the non-null
points. The non-commutative group $G^{r}$ naturally acts on the red non-null
projective points, as well as on the red null projective points. Both actions
are transitive.

The group structure of the red rotations can be transferred to give a
\textbf{multiplication} of red non-null projective points, by the rule
\[
\left[  a:b\right]  \times_{r}\left[  c:d\right]  =\left[  ac+bd:ad+bc\right]
.
\]
This multiplication is associative, commutative, has identity $\left[
1:0\right]  ,$ and the inverse of $\left[  a:b\right]  $ is $\left[
a:-b\right]  $.

\subsection{Green isometries}

\begin{theorem}
An isometry of the green p-quadrance has one of the two forms%
\[%
\begin{tabular}
[c]{lllll}%
$\sigma_{\left[  a:b\right]  }^{g}=%
\begin{bmatrix}
0 & a\\
b & 0
\end{bmatrix}
$ &  & \textrm{or} &  & $\rho_{\left[  a:b\right]  }^{g}=%
\begin{bmatrix}
a & 0\\
0 & b
\end{bmatrix}
$%
\end{tabular}
\]
for some green non-null projective point $\left[  a:b\right]  $.
\end{theorem}

\begin{proof}
Suppose that $\sigma$ sends $j_{1}\equiv\left[  1:1\right]  $ to $\left[
a:b\right]  $ and $j_{2}\equiv\left[  1:-1\right]  $ to $\left[  c:d\right]
$. Then since $q_{g}\left(  j_{1},j_{2}\right)  =1,$ we must have%
\[
-\frac{\left(  ad-bc\right)  ^{2}}{4abcd}=1
\]
which implies that
\[
ad+bc=0.
\]
So $\left[  c:d\right]  =\left[  a:-b\right]  .$ Now given an arbitrary
p-point $u\equiv\left[  x:y\right]  $ with $u\sigma\equiv v\equiv\left[
w:z\right]  $,%
\begin{align*}
q_{g}\left(  u,j_{1}\right)   & =-\frac{\left(  x-y\right)  ^{2}}{4xy}\\
& =q_{g}\left(  v,\left[  a:b\right]  \right)  =-\frac{\left(  za-wb\right)
^{2}}{4wzab}%
\end{align*}
and%
\begin{align*}
q_{g}\left(  u,j_{2}\right)   & =\frac{\left(  x+y\right)  ^{2}}{4xy}\\
& =q_{g}\left(  v,\left[  a:-b\right]  \right)  =\frac{\left(  wb+za\right)
^{2}}{4wzab}.
\end{align*}
So
\begin{align*}
-\frac{\left(  x-y\right)  ^{2}}{4xy}  & =-\frac{\left(  za-wb\right)  ^{2}%
}{4wzab}\\
\frac{\left(  x+y\right)  ^{2}}{4xy}  & =\frac{\left(  wb+za\right)  ^{2}%
}{4wzab}%
\end{align*}
The two possible solutions for $\left[  w:z\right]  $ are%
\[%
\begin{tabular}
[c]{lllll}%
$\left[  ay:bx\right]  $ &  & \textrm{and} &  & $\left[  ax:by\right]  $%
\end{tabular}
\]
which correspond to $u\sigma_{\left[  a:b\right]  }^{g}$ and $u\rho_{\left[
a:b\right]  }^{g}$ respectively.
\end{proof}

\begin{theorem}
For any non-null projective points $\left[  a:b\right]  $ and $\left[
c:d\right]  $%
\[%
\begin{tabular}
[c]{lll}%
$\sigma_{\left[  a:b\right]  }^{g}\sigma_{\left[  c:d\right]  }^{g}%
=\rho_{\left[  ad:bc\right]  }^{g}$ &  & $\rho_{\left[  a:b\right]  }^{g}%
\rho_{\left[  c:d\right]  }^{g}=\rho_{\left[  ac:bd\right]  }^{g}$\\
&  & \\
$\rho_{\left[  a:b\right]  }^{g}\sigma_{\left[  c:d\right]  }^{g}%
=\sigma_{\left[  ac:bd\right]  }^{g}$ &  & $\sigma_{\left[  a:b\right]  }%
^{g}\rho_{\left[  c:d\right]  }^{g}=\sigma_{\left[  ad:bc\right]  }^{g}.$%
\end{tabular}
\]

\end{theorem}

\begin{proof}
These are immediate. Note that the resultant isometries are also associated to
non-null points.
\end{proof}

We now define the group $G^{g}$ of all isometries of the green quadrance, and
distinguish the subgroup $G_{e}^{g}$ of green rotations $\rho_{\left[
a:b\right]  }^{g}.$ These latter are naturally in bijection with the non-null
projective points.

The coset $G_{o}^{g}$ consists of green reflections $\sigma_{\left[
a:b\right]  }^{g}$ which are also naturally in bijection with the non-null
points. The group $G^{g}$ naturally acts on the green non-null projective
points, as well as on the green null projective points. Both actions are transitive.

The group structure of the green rotations can be transferred to give a
\textbf{multiplication} of green non-null projective points, by the rule
\[
\left[  a:b\right]  \times_{g}\left[  c:d\right]  =\left[  ac:bd\right]  .
\]
This multiplication is associative, commutative, has identity $\left[
1:1\right]  ,$ and the inverse of $\left[  a:b\right]  $ is $\left[
b:a\right]  $. This is a familiar algebraic object: it is just the
multiplicative group of non-zero fractions.

In particular if we wish to take powers in this group, then say $\left[
2:3\right]  ^{2}=\left[  4:9\right]  $ and $\left[  2:3\right]  ^{3}=\left[
8:27\right]  $ and so on. But each of these powers are then equally spaced
with respect to green p-quadrance. As an application we may easily prove the
following result about values of spread polynomials.

\begin{theorem}
If $s=-\frac{\left(  y-x\right)  ^{2}}{4xy}$ then $S_{n}\left(  s\right)
=-\frac{\left(  y^{n}-x^{n}\right)  ^{2}}{4x^{n}y^{n}}.$
\end{theorem}

\begin{proof}
The number $s=-\frac{\left(  y-x\right)  ^{2}}{4xy}$ is $q_{g}\left(
1,1,x,y\right)  $ while $r=-\frac{\left(  y^{n}-x^{n}\right)  ^{2}}%
{4x^{n}y^{n}} $ is $q_{g}\left(  1,1,x^{n},y^{n}\right)  $. But since in the
green multiplication $\left[  x:y\right]  ^{n}=\left[  x^{n}:y^{n}\right]  $
we must have $r=S_{n}\left(  s\right)  .$
\end{proof}

\end{document}